\documentclass[a4paper,10pt]{article}
\usepackage{graphicx}
\usepackage{latexsym}
\usepackage{amsmath}
\usepackage{amsthm}
\usepackage{amssymb}
\usepackage{color}
\newcommand{\ol}[1]{\overline{#1}}
\newtheorem{thm}{Theorem}
\newtheorem{cor}{Corollary}
\newtheorem{claim}{Claim}
\newtheorem{lm}{Lemma}

\def\qed{\hfill$\Box$}

\def\P{\mathcal{P}}
\def\D{\mathcal{D}}

\allowdisplaybreaks
\begin{document}
\vspace*{-0.5mm}
\begin{center}
{\LARGE{}Distance restricted matching extensions
\medskip\\
in regular non-bipartite graphs}
\vspace*{0.5cm}\\
{\large
\begin{center}
Jun Fujisawa\footnotemark\\
Faculty of Business and Commerce,
Keio University\\
Hiyoshi 4--1--1, Kohoku-Ku,
Yokohama, Kanagawa 223--8521,
Japan\\
\texttt{fujisawa6@keio.jp}
\end{center}
\footnotetext{{}work supported by JSPS KAKENHI Grant Numbers JP20K03723 and JP24K06833}}
\end{center}
\medskip
\begin{abstract}
Let $m$ and $r$ be integers with $m \ge r \ge 3$ and
let $G$ be an $r$-regular graph of even order.
Let $M$ be a matching in $G$ of size $m$ such that
each pair of edges in $M$ is at distance at least $3$.
In 2023, Aldred et al.~proved that
if $G$ is cyclically $(mr-r+1)$-edge-connected and $G$ is bipartite,
then there exists a perfect matching of $G$ containing $M$.
In this paper, we present non-bipartite analogues of Aldred et al.'s theorem.
An odd ear of $U \subseteq V(G)$ is a path of odd length 
whose ends lie in $U$ but whose internal vertices do not,
or a cycle of odd length having exactly one vertex in $U$.
Our first result shows that
if $G$ is cyclically $(mr - m +1)$-edge-connected and
there exist $mr - \left\lceil \frac{r}{2} \right\rceil + 1$
edge-disjoint odd ears of $V(M)$,
then $M$ can be extended to a perfect matching of $G$.
We further show that if $G$ contains $mr-r+1$ edge-disjoint odd ears of $V(M)$ and
no cyclic edge cut in $G$ of size less than
$(2m-1)(r-1)$ separates an odd cycle from another cycle,
then $M$ can still be extended to a perfect matching.
The second result extends Aldred et al.'s theorem
to non-bipartite graphs in the case $r \ge 4$,
and in the case when $r = 3$ and each pair of edges in $M$ is at distance at least $5$.
It is also shown that
the above results hold when $m \le r - 1$,
without assuming the distance condition on $M$.
\end{abstract}
\vspace{0.5\baselineskip}
\section{Introduction}
\label{intro}
Throughout this paper we consider only finite simple graphs.
A set $F$ of edges in a graph $G$ is called a \textit{cyclic edge cut}
if $G-F$ has two components each of which contains a cycle.
The \textit{cyclic edge-connectivity} of $G$, denoted by $\lambda_c(G)$,
is defined as the size of a minimum cyclic edge cut of $G$.
If there exists no cyclic edge cut in $G$, then we let $\lambda_c(G)= \infty$.
For edges $e,f \in E(G)$, we define the distance between $e$ and $f$
to be the length of the shortest path joining an endvertex of $e$ and an endvertex of $f$.
\medskip

The study of perfect matchings in cubic graphs has received much attention in graph theory.
One of the problems that attracted considerable interest
is the Lov\'{a}sz-Plummer Conjecture, proposed in the 1970s,
stating that the number of perfect matchings in a bridgeless cubic graph
grows exponentially with its order (see \cite{LP}).
The conjecture was proven for bipartite graphs by Voorhoeve \cite{Vo},
for planar graphs by Chudnovsky and Seymour \cite{CS},
and the full conjecture was proven by Esperet et al.~\cite{EKKKN}.
Another well-known conjecture, widely referred to as the Berge-Fulkerson Conjecture,
states that every bridgeless cubic graph $G$ admits six perfect matchings
such that each edge in $G$ is contained in exactly two of them.
Regarding this conjecture, even the following weaker form,
conjectured by Fan and Raspaud,
remains widely open: every bridgeless cubic graph admits three perfect matchings
$M_1$, $M_2$, $M_3$ such that $M_1 \cap M_2 \cap M_3 = \emptyset$.
For recent progress related to this topic, the reader is referred to \cite{KMZ}.
\medskip

In a graph $G$,
we say that a matching $M \subseteq E(G)$ \textit{extends} to a perfect matching 
if there exists a perfect matching $M^*$ of $G$ containing $M$.
Studying matching extension in cubic graphs, and more generally in regular graphs,
is of importance, as it provides explicit information about the structure of their perfect matchings.
A classical result due to K\"{o}nig \cite{Ko1, Ko2} states that
if $G$ is an $r$-regular bipartite graph, then $G$ is $r$-edge-colorable;
hence every edge of $G$ extends to a perfect matching.
However, a matching with two edges in an $r$-regular bipartite graph
does not necessarily extend to a perfect matching (see \cite{AJP}).
Concerning the extendability of large matchings,
Plummer proved the following theorem.

\begin{thm}[Plummer \cite{Pl91}]
Let $m$ and $r$ be positive integers with $m \le r-1$ and
let $G$ be an $r$-regular bipartite graph with $\lambda_c(G) \ge mr-r+1$.
Then every matching $M$ of $G$ with $|M| = m$ extends to a perfect matching.
\label{Pl91}
\end{thm}

If $M$ is a matching of size at least $r$ in an $r$-regular graph $G$
such that $G-V(M)$ has a component consisting of a single vertex,
then $M$ cannot be extended to a perfect matching.
To exclude such cases,
Aldred and Jackson imposed a proximity condition on the matching to be extended,
and investigated the extendability of matchings with three or more edges
in cubic bipartite graphs.
A matching $M$ is said to be a \textit{distance $d$ matching}
if the distance between each pair of distinct edges in $M$ is at least $d$.

\begin{thm}[Aldred and Jackson \cite{AJ}]
Let $G$ be a cubic bipartite graph with $\lambda_c(G) \ge 3m-2$
for some positive integer $m$.
If $M$ is a distance $5$ matching in $G$ with $|M| = m$,
then $M$ extends to a perfect matching.
\label{AJ}
\end{thm}

Theorem \ref{AJ} is generalized to $r$-regular bipartite graphs
for all $r \ge 3$, as follows.

\begin{thm}[Aldred, Jackson and Plummer \cite{AJP}]
Let $m$ and $r$ be integers with $m \ge r \ge 3$ and
let $G$ be an $r$-regular bipartite graph with $\lambda_c(G) \ge mr-r+1$.
If $M$ is a distance $3$ matching in $G$ with $|M| = m$,
then $M$ extends to a perfect matching.
\label{AJP}
\end{thm}

Let $H$ be an $r$-regular bipartite graph and
let $x_1 y_1$, $x_2 y_2$ be two edges of $H$ such that
$x_1$ and $x_2$ belong to the same partite set of $H$.
Let $G = H - x_1 y_1 - x_2 y_2 + x_1 x_2 + y_1 y_2$,
then $G$ is also $r$-regular (and almost bipartite),
but there exist two edges
$x_1y_1'$ and $y_1y_2$ that cannot be extended to a perfect matching,
where $y_1' \in N_G(x_1) \setminus \{x_2\}$.
Moreover, if the cyclic edge-connectivity of the original graph $H$ is sufficiently large,
then the edges $x_1 y_1$ and $x_2 y_2$ can be chosen so that
the distance between $x_1y_1'$ and $y_1y_2$ in $G$ is arbitrarily large.
Thus, as noted in \cite{AJ},
the lack of biparticity seems to be a critical obstacle
to achieving the desired perfect matching
in Theorems \ref{Pl91}, \ref{AJ} and \ref{AJP}.
However, many important families of non-bipartite regular graphs, such as snarks, do exist,
and hence investigating the matching extendability of such graphs
is of particular interest.
The purpose of this paper is to provide new insights into this topic.
\medskip

Let $G$ be a graph and let $U$ be a subset of $V(G)$.
An \textit{ear} of $U$ is defined as
either a path having both endvertices but no internal vertices in $U$,
or a cycle having exactly one vertex in $U$.
In contrast to the context of ear decompositions,
here we define ears with respect to a vertex subset,
and hence every edge joining two vertices of $U$ serves as an ear of $U$.
An ear is \textit{odd} if it contains odd number of edges.
Our first result shows that an analogue of Theorem \ref{AJP} holds for non-bipartite graphs,
provided that there exist many edge-disjoint odd ears of $V(M)$.
\begin{thm}
Let $m$ and $r$ be integers with $m \ge r \ge 3$.
Let $G$ be an $r$-regular graph of even order with
$\lambda_c(G) \ge mr - m + 1$,
and let $M$ be a distance $3$ matching in $G$ with $|M| = m$.
If there exist $mr - \left\lceil \frac{r}2 \right\rceil + 1$
edge-disjoint odd ears of $V(M)$ in $G$, then $M$ extends to a perfect matching.
\label{main-ear}
\end{thm}

Note that, in the case where $m > r$,
the condition on $\lambda_c(G)$ in Theorem \ref{main-ear}
is weaker than that in Theorems \ref{AJ} and \ref{AJP}.
\medskip

The following lemma shows that, in regular bipartite graphs,
large cyclic edge-connectivity guarantees
the existence of many edge-disjoint odd ears of $V(M)$,
except for some exceptional cases.
\begin{lm}
Let $d$, $m$ and $r$ be integers with $d$, $r \ge 3$,
let $G$ be an $r$-regular bipartite graph and 
let $M$ be a distance $d$ matching in $G$ with $|M| = m$.
If $(r,d) \neq (3,3), (3,4)$,
then there exist $\min\{\lambda_c(G), mr\}$ edge-disjoint odd ears of $V(M)$ in $G$.
\label{ears}
\end{lm}

Theorem \ref{main-ear} and Lemma \ref{ears} yield the following corollary.
\begin{cor}
Let $d$, $m$ and $r$ be integers with $d \ge 3$, $m \ge r \ge 3$
and $(r,d) \neq (3,3), (3,4)$.
Let $G$ be an $r$-regular bipartite graph with
$\lambda_c(G) \ge mr - \left\lceil \frac{r}2 \right\rceil + 1$.
If $M$ is a distance $d$ matching in $G$ with $|M| = m$,
then $M$ extends to a perfect matching in $G$.
\label{cor-ear}
\end{cor}
The condition on $\lambda_c(G)$ in Corollary \ref{cor-ear} differs by
$\lfloor \frac{r}2 \rfloor$ from that in Theorems \ref{AJ} and \ref{AJP}.
To fill this gap, we introduce an additional condition.
We call $F$ an \textit{odd-cyclic edge cut}
if $F$ is a cyclic edge cut such that a component of $G-F$ contains a cycle of odd length.
The size of a minimum odd-cyclic edge cut
is denoted by $\lambda_{oc}(G)$.
If there exists no odd-cyclic edge cut in $G$,
then we let $\lambda_{oc}(G)= \infty$.
Our second result shows that,
by imposing an assumption on $\lambda_{oc}(G)$ instead of $\lambda_c(G)$,
the assumption on odd ears in Theorem \ref{main-ear} can be relaxed.
\begin{thm}
Let $m$ and $r$ be integers with $m \ge r \ge 3$.
Let $G$ be an $r$-regular graph of even order with
$\lambda_{oc}(G) \ge (2m-1)(r-1)$,
and let $M$ be a distance $3$ matching in $G$ with $|M| = m$.
If there exist $mr-r+1$ edge-disjoint odd ears of $V(M)$ in $G$,
then $M$ extends to a perfect matching.
\label{main-oc}
\end{thm}

Note that $\lambda_{oc}(G)= \infty$ for any bipartite graph $G$.
Therefore, it follows from Lemma \ref{ears} that
Theorem \ref{main-oc} generalizes Theorem \ref{AJ} to non-bipartite graphs.
Moreover, Theorem \ref{main-oc} generalizes Theorem \ref{AJP},
with the exception of the case where $r=3$ and $M$ is not a distance $5$ matching.
\medskip

Furthermore, we show that
Theorems \ref{main-ear} and \ref{main-oc} also hold in the case $m \le r - 1$,
without assuming the distance condition on $M$.
\begin{thm}
Let $m$ and $r$ be integers with $2 \le m \le r-1$.
Let $G$ be an $r$-regular graph of even order with
$\lambda_c(G) \ge mr - m + 1$,
and let $M$ be a matching in $G$ with $|M| = m$.
If there exist $mr - \left\lceil \frac{r}2 \right\rceil + 1$
edge-disjoint odd ears of $V(M)$ in $G$,
then $M$ extends to a perfect matching.
\label{main2-ear}
\end{thm}
\begin{thm}
Let $m$ and $r$ be positive integers with $2 \le m \le r-1$.
Let $G$ be an $r$-regular graph of even order with
$\lambda_{oc}(G) \ge (2m-1)(r-1)$,
and let $M$ be a matching in $G$ with $|M| = m$.
If there exist $mr-r+1$ edge-disjoint odd ears of $V(M)$ in $G$,
then $M$ extends to a perfect matching.
\label{main2-oc}
\end{thm}

Unlike Theorem 1, Theorem 6 does not hold under the condition
$\lambda_c(G) \ge mr-r+1$.
By a similar argument as above, it follows that
Theorem 7 generalizes Theorem 1 to non-bipartite graphs,
with the exception of the case where $r=3$ and $M$ is not a distance $5$ matching.
\medskip

The paper is organized as follows. 
In the rest of this section
we introduce some terminology and notation.
In Section \ref{Proof},
we prove Theorems \ref{main-ear}--\ref{main2-oc} in a unified form,
under slightly weaker conditions.
In Section \ref{sharpness},
we demonstrate the sharpness of this theorem.
Lemma \ref{ears} is proved in Section \ref{relation}.
\medskip

Throughout this paper,
we use $\ol{H}$ to denote $G-V(H)$
for any subgraph $H$ of a graph $G$.
For $X, Y \subseteq V(G)$,
the set of edges joining $X$ and $Y$ is denoted by $E(X,Y)$.
We often identify a graph with its vertex set when there is no fear of confusion.
For example, if $H_1$ and $H_2$ are subgraphs of $G$,
we use $E(H_1,H_2)$ instead of $E(V(H_1), V(H_2))$.
A path $P = v_1 v_2 \ldots v_k$ is called an $(X,Y)$-path  
if $ V(P) \cap X = \{v_1\} $ and $ V(P) \cap Y = \{v_k\}$.
\medskip

For an integer $d$,
we call $U \subseteq V(G)$ a \textit{distance $d$ set}
if each pair of vertices in $U$ is at distance at least $d$ in $G$.
Let $\ell$, $d$ and $r$ be positive integers with $r \ge 3$.
It follows from Theorem 6 of \cite{LH} that
there exists an $r$-regular bipartite graph $G(\ell,d;r)$ with $\lambda_c(G(\ell,d;r)) \ge \ell(d+2)(r-2)$.
Note that the girth of $G(\ell,d;r)$ is at least $\ell(d+2)$.
Hence, around a shortest cycle of $G(\ell,d;r)$,
we can take a distance $d+1$ matching $M$ of size $\ell$.
We denote by $\Gamma(\ell,d;r)$ the graph obtained from $G(\ell,d;r)$ by removing the edges of $M$.
Note that $\Gamma(\ell,d;r)$ has $2\ell$ vertices of degree $r-1$ and
all the other vertices have degree $r$.
Moreover, the set of vertices of degree $r-1$ is a distance $d+1$ set,
and $\lambda_c(\Gamma(\ell,d;r)) \ge \ell(d+2)(r-2) - \ell > d$.
We use this graph in Sections \ref{relation} and \ref{sharpness}.
\section{Proof of Theorems \ref{main-ear}, \ref{main-oc}, \ref{main2-ear} and \ref{main2-oc}}
\label{Proof}
In \cite{AJP}, Theorem \ref{AJP} is in fact established
with a weaker distance constraint on the edges of $M$.
Theorems \ref{main-ear} and \ref{main-oc} can similarly be refined.
Moreover, the assumption on $\lambda_c(G)$
in Theorems \ref{main-ear} and \ref{main2-ear} is slightly relaxed
when either $m$ or $r$ is odd.
With these improvements applied, we present Theorems \ref{main-ear}--\ref{main2-oc}
in a unified statement:
\begin{thm}
Let $k$, $m$ and $r$ be integers with $m \ge 2$ and $r \ge 3$.
Let $G$ be an $r$-regular graph of even order,
and let $M$ be a matching in $G$ with $|M| = m$ such that
there exist $k$ edge-disjoint odd ears of $V(M)$.
If $m \ge r$, we further assume that
no vertex of $V(G) \setminus V(M)$ is adjacent to $r-1$ vertices of $V(M)$.
If one of the following holds, then $M$ extends to a perfect matching:
\begin{itemize}
\item[\textup{(i)}]
$k = mr - \left\lceil \frac{r}2 \right\rceil + 1$ and
$\lambda_c(G) \ge mr-m+\theta$,
where $\theta = 1$ when both $m$ and $r$ are even
and $\theta = 0$ otherwise, or
\item[\textup{(ii)}]
$k = mr-r+1$ and $\lambda_{oc}(G) \ge (2m-1)(r-1)$.
\end{itemize}
\vspace*{-2mm}
\label{main}
\end{thm}
\noindent
\textit{Proof.}\quad
Suppose to the contrary that
$M$ is a matching of $G$ which satisfies the assumptions of the theorem
and does not extend to a perfect matching.
Then by Tutte's theorem there is a barrier set $S \subseteq V(G - V(M))$
such that $G - V(M)-S$ contains at least $|S|+2$ odd components.
Let $T=G-V(M)-S$ and
let $\D$ be the set of components of $T$.
\begin{claim}
If $D$ is an odd component of $T$, then $|E(D, \ol{D})| \ge r$.
\label{t-sanhen}
\end{claim}
\noindent
\textit{Proof.}\quad
Suppose to the contrary that $|E(D, \ol{D})| < r$.
Since $G$ is $r$-regular, 
$|E(H, \ol{H})| = r|V(H)| - 2|E(H)| \ge r$
for every induced forest $H$ in $G$.
Hence neither $D$ nor $\ol{D}$ is a forest,
that is, $E(D, \ol{D})$ is a cyclic edge-cut.
If (i) of the statement of the theorem holds,
then since $m \ge 2$ and $r \ge 3$, we obtain
$|E(D, \ol{D})| \ge \lambda_c(G) \ge mr-m + \theta  > r$, a contradiction.
Hence (ii) of the statement of the theorem holds.
Since $\lambda_{oc}(G) \ge  (2m-1)(r-1) > r$, $D$ does not contain an odd cycle,
i.e.~$D$ is bipartite.
Let $B,W$ be the partite sets of $D$.
Since $D$ is an odd component,
we may assume without loss of generality that $|B| > |W|$.
Then
$|E(D, \ol{D})| \ge |E(B, \ol{D})| = \sum_{v \in B} d_G(v)-|E(B, W)|
 \ge \sum_{v \in B} d_G(v) - \sum_{v \in W} d_G(v)
=r(|B|-|W|) \ge r$,
a contradiction.
\qed
\medskip
\begin{claim}
If $D$ is a component of $T$ which contains a cycle,
then there exists a cycle in $G - V(D)$.\label{twoedges}\end{claim}\noindent
\textit{Proof.}\quad
Let $U$ be the set of vertices in $T - V(D)$
which have at least $r-1$ neighbors in $V(M)$.
Suppose that $|U| \ge 2$.
Then, since $U \neq \emptyset$, we have $m \le r-1$ by hypothesis of the theorem.
Since the graph induced by $V(M) \cup U$
contains at least $|U|(r-1) + m \ge 2m + |U| = |V(M) \cup U|$ edges,
we obtain the desired cycle.
Therefore, we may assume that $|U| \le 1$.
\medskip

Let $\D_0= \{D_0 \in \D \mid D_0 \neq D \mbox{ and } V(D_0) \cap U = \emptyset\}$,
then we have $|\D_0| = |\D| - |U| - 1$.
Suppose that $\D_0 = \emptyset$.
Then, since $|\D| \ge |S|+2$,
we have $S = \emptyset$, $|U| = 1$ and $|\D| = 2$.
Let $D^*$ be the component of $T$ which is distinct from $D$.
If $|V(D^*)|\ge 2$, then since $G$ is $r$-regular and $S = \emptyset$,
every vertex in $V(D^*) \setminus U$ has at least two neighbors in $D^*$,
and thus we obtain the desired cycle in $D^*$.
On the other hand, if $V(D^*)$ consists of one vertex, say $u$,
then since $S = \emptyset$, $u$ is adjacent to $r$ vertices of $V(M)$.
Since $U \neq \emptyset$, we have $m \le r-1$ by hypothesis,
and hence we obtain the desired cycle in the graph induced by $\{u\} \cup V(M)$.
Therefore, we may assume that $\D_0 \neq \emptyset$.
\medskip

Let $D_0 \in \D_0$.
Since any cycle of $D_0$ is a desired cycle,
we may assume that $D_0$ is a tree.
Then either $D_0$ is trivial or $D_0$ contains at least two leaves.
Since $V(D_0) \cap U = \emptyset$, we obtain $|E(D_0,S)| \ge 2$ in either case.
Let $X = \bigcup_{D_0 \in \D_0} V(D_0)$.
Since $|\D_0| = |\D| - |U| - 1 \ge |S|$,
we have
$\left| E\left(X , S \right) \right| \ge 2|\D_0| \ge |\D_0| + |S|$,
and hence there exists a desired cycle in the graph induced by $X \cup S$.
\qed
\medskip

Let $|S| = s$, $|E(S,S)|+|E(S,V(M))| = \mu$ and $|E(V(M),V(M))| = m + m^*$.
Then
\begin{eqnarray}
|E(T, \ol{T})| & = & |E(S, T)| + |E(V(M), T)|\nonumber\\
& = & sr - |E(S,V(M))| - 2|E(S,S)| + 2mr - |E(V(M),S)| - 2|E(V(M),V(M))|\nonumber\\
& = & sr + 2(mr - m - m^* - \mu).
\label{t_tbar}
\end{eqnarray}

Let $P_1, \ldots , P_k$ be edge-disjoint odd ears of $V(M)$.
Without loss of generality, we may assume that each of
$P_1, \ldots , P_{k-m-m^*-\mu}$ does not contain
any edge of $E(V(M),V(M)) \cup E(S,V(M)) \cup E(S,S)$.
Then, since $|E(P_j)|$ is odd,
$|E(P_j) \cap E(T,T)|$ is odd
for every $j$ with $1 \le j \le k-m-m^*-\mu$.
Thus we can take an edge $f_j \in E(P_j) \cap E(T,T)$
so that the component of $P_j \cap T$ which contains $f_j$
is a path of odd length.
Let $F = \{f_j \mid 1 \le j \le k-m-m^*-\mu\}$.
\medskip

\begin{claim}
For every component $D$ of $T$,
$|E(D, \ol{D})| \ge 2|E(D) \cap F|$.
Moreover, if $D$ is a bipartite odd component, then
$|E(D, \ol{D})| \ge 2|E(D) \cap F| + r$.
\label{d-dbar}
\end{claim}
\noindent
\textit{Proof.}\quad
If $f_j \in E(D)$ for some $j$,
then $|E(D, \ol{D}) \cap E(P_j)| \ge 2$.
Since $P_1, \ldots , P_k$ are edge-disjoint,
we have $|E(D, \ol{D})| \ge 2|E(D) \cap F|$.
\medskip

Suppose that $D$ is a bipartite odd component.
Let $B,W$ be the partite sets of $D$, where $|B| > |W|$.
If $f_j \in E(D)$ for some $j$, then
we obtain $|E(P_j) \cap E(W, \ol{D})| \ge 1$,
since the component of $P_j \cap T$ which contains $f_j$
is a path of odd length.
Hence $|E(W, \ol{D})| \ge |E(D) \cap F|$.
Since $|E(W, \ol{D})|+|E(W, B)| = \sum_{v \in W} d_G(v) = r|W|$ and
$|E(B, \ol{D})|+|E(B, W)| = \sum_{v \in B} d_G(v) = r|B| \ge r(|W|+1)$,
we obtain
$|E(B, \ol{D})| \ge |E(W, \ol{D})| + r \ge |E(D) \cap F|+r$.
Hence
$|E(D, \ol{D})| = |E(W, \ol{D})|+|E(B, \ol{D})|
\ge 2|E(D) \cap F| + r$.
\qed
\medskip
\begin{claim}
At most one odd component of $T$ is non-bipartite.
\label{oddcycle}
\end{claim}
\noindent
\textit{Proof.}\quad
Suppose, to the contrary, that there exist two non-bipartite components,
say $D_1$ and $D_2$, of $T$.
Then both $D_1$ and $D_2$ contain an odd cycle.
If (i) of the statement of the theorem holds,
then $|E(D_i, \ol{D_i})| \ge \lambda_c(G)$ for $i=1,2$.
Moreover, since $|V(D_i)|$ is odd, we have
$r \equiv |E(D_i, \ol{D_i})| \pmod 2$,
and hence $|E(D_i, \ol{D_i})| \neq m(r-1)$
if either $m$ or $r$ is odd.
Thus we obtain
$|E(D_i, \ol{D_i})| \ge mr - m + 1$
regardless of the parity of $m$ and $r$.
On the other hand, if (ii) of the statement of the theorem holds,
then $|E(D_i, \ol{D_i})| \ge \lambda_{oc}(G) \ge (2m-1)(r-1) > mr - m + 1$ for $i=1,2$.
In either case, it follows from Claim \ref{t-sanhen} that
$|E(T, \ol{T})| = \sum_{D \in \D} |E(D, \ol{D})|
\ge 2(mr - m + 1) + sr > 2mr - 2m + sr$,
which contradicts (\ref{t_tbar}).
\qed
\medskip

Let $q_1$ and $q_2$ be the number of bipartite odd components
and non-bipartite odd components of $T$, respectively.
Then we have $q_1 + q_2 \ge s+2$,
and it follows from Claim \ref{oddcycle} that $q_2 \le 1$.
By Claim \ref{d-dbar}, we obtain
\begin{eqnarray}
|E(T, \ol{T})|
 & \ge & \sum_{D \in \D} 2|E(D) \cap F| + q_1 r \nonumber\\
& = & 2(k - m - m^* - \mu)  + q_1 r \label{girifour}\\
&\ge& 2(k - m - m^* - \mu) + (s + 1)r. \nonumber
\end{eqnarray}
If (i) of the statement of the theorem holds,
then we have $k > mr - \frac{r}2$.
Hence $|E(T, \ol{T})| > 2(mr - m - m^* - \mu) + sr$,
which contradicts (\ref{t_tbar}).
Thus (ii) of the statement of the theorem holds,
that is, $k = mr-r + 1$ and
$\lambda_{oc}(G) \ge (2m - 1)(r-1)$.
\medskip

If $q_1 \ge s+2$, then (\ref{girifour}) yields
$|E(T, \ol{T})| \ge 2(mr - m - m^* - \mu) + sr + 2$,
which contradicts (\ref{t_tbar}).
Hence we have $q_1 = s+1$ and $q_2 = 1$.
Let $D_1$ be the unique non-bipartite odd component of $T$.
By Claim \ref{twoedges},
there exists a cycle in the graph induced by $G-V(D_1)$,
and hence $|E(D_1, \ol{D_1})| \ge \lambda_{oc}(G) \ge (2m - 1)(r-1)$.
Then it follows from Claim \ref{d-dbar} that
$|E(T, \ol{T})| = \sum_{D \in \D} |E(D, \ol{D})|
 \ge q_1 r + (2m - 1)(r-1) = sr + 2m(r-1) + 1$,
which contradicts (\ref{t_tbar}).
This completes the proof of Theorem \ref{main}.
\qed
\medskip

\section{Relationship between the new and previous results}
\label{relation}

In this section we prove Lemma \ref{ears},
which reveals the relationship between our new theorems and
Theorems \ref{Pl91}--\ref{AJP}.
A distance $d$ set $U$ in a graph $G$ is called \textit{good}
if every path joining two vertices of $U$
contains at least $d-1$ internal vertices of degree at least $3$ in $G$.
Note that, in the graph $\Gamma(l,d;r)$ defined at the end of Section \ref{intro},
the set of vertices of degree $r-1$ is a good distance $d$ set.
\medskip

The following lemma plays a central role in the proof of Lemma \ref{ears}.
\begin{lm}
Let $r$ and $d$ be positive integers with $(r-2)(d-2)\ge 4$.
Let $G$ be a connected $r$-regular graph, let $T$ be an induced tree of $G$
and let $L$ be a subset of $V(T)$.
If $L$ is a distance $d$ set in $G$,
then $|E(T, \ol{T})| \ge r|L|$.
\label{edgecount}
\end{lm}
\noindent
\textit{Proof.}\quad
Suppose to the contrary that there exists a counterexample $G$.
Then, for an induced tree $T$ of $G$,
there exists $L \subseteq V(T)$
such that $L$ is a distance $d$ set in $G$ and $|E(T, \ol{T})| < r|L|$.
Take such $T$ so that $|V(T)|$ is minimum,
and take $L = \{b_1, \ldots , b_{\ell}\} \subseteq V(T)$ as above.
Note that $(r-2)(d-2)\ge 4$ implies $r \ge 3$.
\medskip

Since $G$ is $r$-regular,
$|E(H,\ol{H})|= r|V(H)|-2|E(H)|= (r-2)|V(H)|+2$
for every induced tree $H$ in $G$.
Hence
\[(r-2)|V(T)|+2 = |E(T, \ol{T})| < r|L| = r{\ell},\]
which yields ${\ell} \neq 1$.
If ${\ell} = 2$, then we have $|V(T)| \ge d+1$
since the distance between $b_1$ and $b_2$ is at least $d$.
However, $(r-2)(d-2)\ge 4$ yields
$(r-2)|V(T)|+2 \ge (r-2)(d+1)+2 = (r-2)(d-2) + 3r - 4 > 2r = r{\ell}$, a contradiction.
Hence we have ${\ell} \ge 3$.
\medskip

Suppose that $b_i$ is not a leaf of $T$ for some $i$.
Let $a_1, \ldots , a_k$ be the neighbors of $b_i$ in $T$,
let $T_1$ be the component of $T-\{a_2, \ldots , a_k\}$ which contains $b_i$
and let $T_2$ be the component of $T-a_1$ which contains $b_i$.
Then $V(T_1) \cap V(T_2) = \{b_i\}$.
By the minimality of $T$, we have
$(r-2)|V(T_j)|+2 = |E(T_j, \ol{T_j})| \ge r|V(T_j) \cap L|$ for $j=1,2$.
Since $|V(T_1) \cap L| + |V(T_2) \cap L| = |V(T) \cap L| + 1$,
we obtain
\begin{eqnarray*}
(r-2)|V(T)|+2
& = & (r-2)(|V(T_1)| +|V(T_2)| -1)+2\\
& = & (r-2)|V(T_1)|+2 +(r-2)|V(T_2)|+2 -r \\
& \ge & r|V(T_1) \cap L| + r|V(T_2) \cap L| -r
 \ = \ r|V(T) \cap L|
 \ = \  r{\ell},
\end{eqnarray*}
a contradiction.
Thus every vertex of $L$ is a leaf of $T$.
\medskip

Let us consider $T$ as a rooted tree with the root $b_1$.
For every vertex $v$ of $T$,
let $T_v$ be the subtree of $T$ which is induced by $v$ and its descendants,
and let $L(v) = |V(T_v) \cap L|$.
Moreover, let $x$ be the vertex of $T$ with $L(x) \ge 2$
that is farthest from $b_1$ in $T$.
Since $L(x) \ge 2$, $x$ is not a leaf of $T$, and hence $x \notin L$.
Thus $L(x)=\sum_{x^+ \in N^+(x)} L(x^+)$,
where $N^+(x)$ is the set of children of $x$.
By the choice of $x$, we have
$L(x^+) \le 1$ for every $x^+ \in N^+(x)$,
and hence there exist $u, w \in N^+(x)$ such that $L(u) = L(w) = 1$.
Since each of $T_{u}$ and $T_{w}$ contains a vertex of $L$ and
$L$ is a distance $d$ set,
we have $|V(T_{u})|+|V(T_{w})| \ge d$.
Let $T' = T-V(T_{u})-V(T_{w})$.
Then by the minimality of $T$,
$(r-2)|V(T')|+2 = |E(T', \ol{T'})| \ge r |V(T') \cap L| = r({\ell}-2)$.
Hence
\begin{eqnarray*}
(r-2)|V(T)|+2
& = & (r-2)(|V(T_{u})| + |V(T_{w})|) + (r-2)|V(T')|+2\\
& \ge & (r-2)d + r({\ell} - 2)\\
& = & r{\ell} + (r-2)(d-2) - 4 \ \ge \ r{\ell},
\end{eqnarray*}
a contradiction.
\qed
\medskip

\noindent
\textit{Proof of Lemma \ref{ears}.}\quad
Let $(B, W)$ be the bipartition of $G$ and
let $M = \{b_1 w_1 , \ldots , b_m w_m\}$
such that $b_i \in B$ and $w_i \in W$ for every $i$.
Moreover, let $M_B=\{b_1, \ldots , b_m\}$ and $M_W=\{w_1, \ldots , w_m\}$.
Since $G$ is bipartite, every $(M_B, M_W)$-path is of odd length.
Hence it suffices to show that there exist $k$ edge-disjoint $(M_B, M_W)$-paths,
where $k= \min\{\lambda_c(G), mr\}$.
Suppose to the contrary that such paths do not exist.
Then by Menger's theorem, there exists an edge cut $F$ of $G$
with $|F| < k$ which separates $M_B$ and $M_W$.
Since $\lambda_c(G) \ge k$, 
all but at most one component of $G-F$ is a tree.
Without loss of generality, we may assume that
every component of $G-F$ that contains a vertex of $M_B$ is a tree.
Let $T_1, \ldots , T_p$ be these components
and let $M_i = \{b_j w_j \mid b_j \in V(T_i)\}$ for each $i$ with $1 \le i \le p$.
Then $E(T_i, \ol{T_i}) \subseteq F$ for each $i$ and
$\bigcup_{i=1}^p M_i = M$.
Let
$d'=d$ if $d$ is even and $d'=d+1$ if $d$ is odd.
Then we have $(r-2)(d'-1) \ge 4$.
Moreover, since $G$ is bipartite and $M_B \subseteq B$,
$M_B$ is a distance $d'$ set in $G$.
Hence it follows from Lemma \ref{edgecount} that $|E(T_i, \ol{T_i})| \ge r|M_i|$ for every $i$.
Thus we obtain
$|F| \ge \sum_{i=1}^p |E(T_i, \ol{T_i})| \ge \sum_{i=1}^p r|M_i| = r|M| = mr \ge k$,
a contradiction.
\qed
\medskip

The condition $(r,d) \neq (3,3), (3,4)$ in Lemma \ref{ears} is necessary.
This is shown by a cubic bipartite graph $G$ constructed as follows.
Let $k$ be an integer with $k \ge 2$
and let $T'$ be a tree with $V(T') = \{x_i, y_i \mid 1 \le i \le 2k+1\}$
and $E(T') = \{x_i y_i \mid 1 \le i \le 2k+1\} \cup \{x_i x_{i+1} \mid 1 \le i \le 2k \}$.
Let $T$ be the tree obtained from $T'$ by
subdividing each edge of $\{x_{2i} y_{2i} \mid 1 \le i \le k\}$ exactly once.
Note that $T$ has $k+2$ vertices of degree two
and the set of leaves of $T$ is $\{y_1, \ldots , y_{2k+1}\}$.
\medskip

Let $H'$ be the graph isomorphic to $\Gamma(4k+2, 6k+3; 3)$,
the graph defined at the end of Section \ref{intro}.
Let $(B, W)$ be the bipartition of $H'$,
and let
$\{b_i \mid 1 \le i \le 4k+2\}$ and 
$\{w_i \mid 1 \le i \le 4k+2\}$ be the set of vertices of degree two
in $B$ and $W$, respectively.
Let $H$ be the graph obtained from $H'$
by adding $k$ vertices $u_1, \ldots , u_k$ and
adding the set of edges $\{u_i b_j \mid 1 \le i \le k,\ k+3i \le j \le k+3i+2 \}$.
Since $\lambda_c(H') \ge 6k+3$ and $\{b_j \mid k+3 \le j \le 4k+2 \}$
is a good distance $6k+3$ set in $H'$, we have $\lambda_c(H) \ge 6k+3$.
Note that the set of vertices of degree two in $H$ is
$\{b_1 , \ldots , b_{k+2}\} \cup \{w_1, \ldots , w_{4k+2}\}$.
\medskip

Let $G$ be the graph obtained from $T \cup H$ by joining
the $k+2$ vertices of degree two in $T$
and $\{b_1 , \ldots , b_{k+2}\}$ by a matching,
and adding the edge set $\{y_i w_{2i-1}, y_i w_{2i} \mid 1 \le i \le 2k+1 \}$.
Then $G$ is a cubic bipartite graph.
Since $\{b_1 , \ldots , b_{k+2}\} \cup \{w_1, \ldots , w_{4k+2}\}$
is a good distance $6k+3$ set in $H$,
we can deduce that $\lambda_c(G) \ge 6k+3$.
Let $(B_G, W_G)$ be the bipartition of $G$ such that $W \subseteq W_G$,
and let $M = \{y_i w_{2i} \mid 1 \le i \le 2k+1 \}$.
Then $M$ is a distance $4$ matching of $G$ with $|M|=2k+1$.
Note that
$V(M) \cap B_G = \{y_1, \ldots, y_{2k+1}\} \subseteq V(T)$ and
$V(M) \cap W_G = \{w_2, w_4, \ldots, w_{4k+2}\}  \subseteq V(H)$.
Since $G$ is bipartite,
each odd ear of $V(M)$ in $G$ contains an edge which joins $T$ and $H$.
Hence the number of edge-disjoint odd ears of $V(M)$ in $G$
is at most $|E(T,H)|=5k+4$,
which is less than $\min\{\lambda_c(G), 3|M|\}=6k+3$.
Therefore, Lemma \ref{ears} does not hold when $(r,d)= (3,3)$ or $(3,4)$.
\section{Sharpness}
\label{sharpness}
In this section, we show that Theorem \ref{main} is best possible in several ways.
\subsection{Sharpness of the condition on $k$ in (i) and
on $\lambda_{oc}(G)$ in (ii)}
\label{renketsudo}
\def\bup{\lceil \frac{r}2 \rceil}
\def\bdown{\lfloor \frac{r}2 \rfloor}
Let $m$ and $r$ be integers with $m \ge 2$ and $r \ge 3$.
We construct an $r$-regular graph $G$ of even order
with $\lambda_c(G) \ge 2m(r-1)-r$
that contains a non-extendable matching $M$ with $m$ edges.
To this end, we prepare three graphs, $H_1$, $H_2$ and $H_3$.
Let $d$ be an integer with $d \ge 2m(r-1)-r$
and let $\alpha = m(r-1)$.
First, we define $H_1$.
Let $A_1$ be the graph isomorphic to $\Gamma(\alpha + \bup, d; r)$,
and let
$\{x_i \mid 1 \le i \le \alpha + \bup\}$ and
$\{y_i \mid 1 \le i \le \alpha + \bup\}$ be the sets of vertices of degree $r-1$
in the two partite sets of $A_1$, respectively.
Since $A_1$ is obtained
by removing $\alpha + \bup$ edges
from a cycle of an $r$-regular bipartite graph,
we may assume that
\begin{itemize}
\item[(A)]
there exist vertex-disjoint paths $P_1, \ldots , P_{\alpha + \bup}$ in $A_1$
such that $P_i$ joins $x_i$ and $y_i$ for every $i$ with $1 \le i \le \alpha + \bup$.
\end{itemize}
Let $H_1$ be the graph obtained from $A_1$
by adding a vertex $q_1$ and
the set of edges $\{q_1 y_i \mid \alpha - \bdown + 1 \le i \le \alpha + \bup\}$.
Next, we define $H_2$.
Let $A_2$ be the graph isomorphic to $\Gamma(\alpha, d; r)$,
and let
$\{z_i \mid 1 \le i \le \alpha\}$ and
$\{u_i \mid 1 \le i \le \alpha\}$ be the sets of vertices of degree $r-1$
in the two partite sets of $A_2$, respectively.
As in the case of $A_1$, we may assume that
\begin{itemize}
\item[(B)]
there exist vertex-disjoint paths $Q_1, \ldots , Q_{\alpha}$ in $A_2$
such that $Q_i$ joins $z_i$ and $u_i$ for every $i$ with $1 \le i \le \alpha$.
\end{itemize}
Let $H_2$ be the graph obtained from $A_2$
by adding a vertex $q_2$ and the set of edges
$\{q_2 z_i \mid \alpha - \bup + 1 \le i \le \alpha\} \cup
\{q_2 u_i \mid \alpha - \bdown + 1 \le i \le \alpha\}$.
Finally, let $H_3$ be the graph with
$V(H_3) = \{b_i , w_i \mid 1 \le i \le m\}$
and $E(H_3) = \{b_i w_i \mid 1 \le i \le m\}$;
that is, $H_3 \simeq mK_2$.
\medskip

Let
$X = \{x_i \mid 1 \le i \le \alpha + \bup\}$,
$Y = \{y_i \mid 1 \le i \le \alpha - \bdown\}$,
$Z = \{z_i \mid 1 \le i \le \alpha - \bup\}$ and
$U = \{u_i \mid 1 \le i \le \alpha - \bdown\}$.
Note that, in $H_1 \cup H_2$,
every vertex in $X \cup Y \cup Z \cup U$ has degree $r-1$
and all the other vertices have degree $r$.
Moreover,
\begin{itemize}
\item[(C)]
$X \cup Y$ (resp.~$Z \cup U$) is a good distance $d+1$ set in $H_1$ (resp.~$H_2$).
\end{itemize}
We construct $G$ by adding edges to $H_1 \cup H_2 \cup H_3$ as follows:
\begin{itemize}
\item
join $Y$ and $U$ by a matching of size $\alpha - \bdown$,
\item
join $b_i$ to
$\lceil \frac{r-1}2 \rceil$ vertices of $X$ and
$\lfloor \frac{r-1}2 \rfloor$ vertices of $Z$
 for every $i$ with $1 \le i \le m-1$,
\item
join $w_i$ to
$\lfloor \frac{r-1}2 \rfloor$ vertices of $X$ and
$\lceil \frac{r-1}2 \rceil$ vertices of $Z$
for every $i$ with $1 \le i \le m-1$,
\item
join $b_m$ to 
$\bup$ vertices of $X$ and
$\bdown -1$ vertices of $Z$ and
\item
join $w_m$ to $r-1$ vertices of $X$.
\end{itemize}
Note that, in the resulting graph,
$|E(H_3,X)|=(m-1)(r-1)+\bup + r-1 = \alpha + \bup = |X|$ and
$|E(H_3,Z)|=(m-1)(r-1)+\bdown -1 = (m-1)(r-1) + r - \bup - 1 = \alpha - \bup = |Z|$.
Hence,
by joining each vertex of $X \cup Z$ to exactly one vertex of $H_3$ in the above procedure,
we obtain an $r$-regular graph $G$ of even order.
Let $M = E(H_3)$.
Then it follows from (C) that $M$ is a distance $d+3$ matching in $G$.
Let $S$ be the partite set of $H_1$ which contains $Y$.
Then $G - V(M) - S$ consists of the union of $|S|+1$ independent vertices and $H_2$,
and hence $M$ does not extend to a perfect matching.
\medskip

In what follows, we show that $\lambda_c(G) \ge 2m(r-1)-r$.
Recall that $\lambda_c(A_i) \ge d$ for $i=1,2$.
Since $H_i$ is obtained from $A_i$ by adding a vertex $q_i$ and joining it to
the vertices of some good distance $d+1$ set,
we can deduce that $\lambda_c(H_i) \ge d \ge 2m(r-1)-r$ for $i=1,2$.
Note that there exist $|Z|$ edge-disjoint $(X,Z)$-paths of length two or three
in the graph induced by $X \cup Z \cup V(H_3) \setminus \{w_m\}$.
Thus, from another perspective,
$G$ may be regarded as being obtained from $H_1 \cup H_2$
by the following operations:
\begin{itemize}
\item[(a)]
join $Y$ and $U$ by a matching,
\item[(b)]
join $Z$ and $\alpha - \bup$ vertices of $X$
by $\alpha - \bup$ vertex-disjoint paths of length two or three,
\item[(c)]
identify some of the new vertices created in (b)
and denote them by $b_1, \ldots, b_m, w_1, \ldots, w_{m-1}$,
\item[(d)]
join $b_m$ to $\bup - (\bdown - 1)$ vertices of $X$ and
\item[(e)]
add the vertex $w_m$ and join it to $r-1$ vertices of $X$ and $b_m$.
\end{itemize}
In the procedures (a) and (b),
$H_1$ and $H_2$ are joined by $|Y| + \alpha - \bup = 2 \alpha - r = 2m(r-1)-r$
vertex-disjoint paths.
Since $d \ge 2m(r-1)-r$, it follows from (C) that
the cyclic edge-connectivity of the graph
obtained through procedures (a) and (b) is $2m(r-1)-r$.
By (C) again, we can deduce that
the cyclic edge-connectivity remains at least $2m(r-1)-r$
in each step of the procedures (c)---(e).
Therefore, $\lambda_c(G) \ge 2m(r-1)-r$.
\medskip

It follows from (A) and (B) that
$H_1$ contains $\alpha - \bdown$ vertex-disjoint $(X,Y)$-paths of odd length
which cover $\alpha - \bdown$ vertices of $X$ and all the vertices of $Y$,
and
$H_2$ contains $\alpha - \bup$ vertex-disjoint $(Z,U)$-paths of odd length
which cover all the vertices of $Z$ and $\alpha - \bup$ vertices of $U$.
Using these paths together with
some edges in $E(H_1, H_2)$, $E(H_2, H_3)$ and $E(H_3, H_1)$,
we obtain $\alpha-\bup$ edge-disjoint odd ears of $V(M)$ in $G - E(M)$.
Since each edge of $M$ corresponds to an odd ear of $V(M)$,
the total number of edge-disjoint odd ears of $V(M)$ is $\alpha-\bup+ m = mr - \bup$.
Since $\lambda_c(G) \ge 2m(r-1)-r \ge mr-m+1$,
it follows that the condition on $k$ in (i) of Theorem \ref{main} is sharp.
Moreover, the condition $\lambda_{oc}(G) \ge (2m-1)(r-1)$ in (ii) is sharp as well,
since $mr - \bup \ge mr - r + 1$
and $\lambda_{oc}(G) \ge \lambda_c(G) \ge 2m(r-1)-r = (2m-1)(r-1)-1$.

\subsection{Sharpness of the condition on $\lambda_c(G)$ in (i)}
Let $m$ and $r$ be integers with $m \ge 2$ and $r \ge 3$,
and let $\rho = \lceil \frac{(m+1)(r-1)}2 \rceil$.
Let $A$ be a graph isomorphic to $\Gamma(\rho, m(r-1); r)$
and let
$\{x_i \mid 1 \le i \le \rho\}$ and
$\{y_i \mid 1 \le i \le \rho\}$
be the sets of vertices of degree $r-1$
in the two partite sets of $A$, respectively.
As in the argument in Section \ref{renketsudo}, we may assume that
\begin{itemize}
\item[(D)]
there exist vertex-disjoint paths $P_1, \ldots , P_\rho$ in $A$ such that
$P_i$ joins $x_i$ and $y_i$ for every $i$ with $1 \le i \le \rho$.
\end{itemize}
Let $H_1$ be the graph obtained from $A$
by adding a vertex $q$ and adding the set of edges 
$\{q x_i \mid \rho - \lceil \frac{r}2 \rceil + 1 \le i \le \rho\} \cup
\{q y_i \mid \rho - \lfloor \frac{r}2 \rfloor + 1 \le i \le \rho\}$.
Here we prepare a copy of $H_1$ and denote it by $H_2$.
Let $U_i$ be the set of vertices of degree $r-1$ in $H_i$ for $i=1,2$.
Then $|U_1|= |U_2|= 2\rho - r$.
Note that $2\rho - r = m(r-1)$ if both $m$ and $r$ are even,
and $2\rho - r = m(r-1)-1$ if either $m$ or $r$ is odd.
\medskip

Let $H_3$ be the graph with $V(H_3) = \{b_i , w_i \mid 1 \le i \le m\}$
and $E(H_3) = \{b_i w_i \mid 1 \le i \le m\}$.
In the case where both $m$ and $r$ are even,
we add edges to $H_1 \cup H_2 \cup H_3$ as follows:
\begin{itemize}
\item[(f)]
join $b_i$ to
$\lceil \frac{r-1}2 \rceil$ vertices of $U_1$ and
$\lfloor \frac{r-1}2 \rfloor$ vertices of $U_2$
 for every $i$ with $1 \le i \le m$ and
\item[(g)]
join $w_i$ to
$\lfloor \frac{r-1}2 \rfloor$ vertices of $U_1$ and
$\lceil \frac{r-1}2 \rceil$ vertices of $U_2$
for every $i$ with $1 \le i \le m$.
\end{itemize}
By joining each vertex of $U_1 \cup U_2$ to exactly one vertex of $H_3$ in the above procedure,
we obtain an $r$-regular graph, say $G$.
Let $M=E(H_3)$.
Then $M$ does not extend to a perfect matching,
since $G - V(M)$ consists of two odd components, $H_1$ and $H_2$.
By (D),
$H_1$ contains $\rho - \frac{r}2$ edge-disjoint paths of odd length
joining a vertex in
$\{x_i \mid 1 \le i \le \rho - \frac{r}2\}$ and a vertex in
$\{y_i \mid 1 \le i \le \rho - \frac{r}2\}$.
By properly combining these paths with edges in $E(H_1, H_3)$,
we obtain $\rho - \frac{r}2 = \frac{m(r-1)}2$ edge-disjoint odd ears of $V(M)$
using only the edges in $E(H_1, H_3) \cup E(H_1)$.
Similarly, there exist $\frac{m(r-1)}2$ edge-disjoint odd ears of $V(M)$
using only the edges in $E(H_2, H_3) \cup E(H_2)$.
Since each edge of $E(H_3)$ corresponds to an odd ear of $V(M)$,
the total number of edge-disjoint odd ears of $V(M)$ is
$mr > mr - \lceil \frac{r}2 \rceil + 1$.
By the similar argument as in Section \ref{renketsudo},
we can deduce that $\lambda_c(G) = |E(H_1, H_3)| = m(r-1)$.
Hence the condition on $\lambda_c(G)$ in (i) of Theorem \ref{main} is sharp
in the case where both $m$ and $r$ are even.
\medskip

Next, we consider the case where either $m$ or $r$ is odd.
In this case, since $|U_1| = |U_2| = m(r-1) - 1$,
adding edges to $H_1 \cup H_2 \cup H_3$ by procedures (f) and (g)
would result in some vertices in $U_1$ and $U_2$ having degree greater than $r$.
Therefore, in this case, we do not take $H_2$ to be a copy of $H_1$,
and instead construct it as follows.
Let $A'$ be a graph isomorphic to $\Gamma(\rho+1, m(r-1); r)$
and let
$\{x'_i \mid 1 \le i \le \rho+1\}$ and
$\{y'_i \mid 1 \le i \le \rho+1\}$
be the sets of vertices of degree $r-1$
in the two partite sets of $A'$, respectively.
Let $H_2$ be the graph obtained from $A'$
by adding a vertex $q'$ and adding the set of edges 
$\{q' x'_i \mid \rho - \lceil \frac{r}2 \rceil + 2 \le i \le \rho + 1\} \cup
\{q' y'_i \mid \rho - \lfloor \frac{r}2 \rfloor + 2 \le i \le \rho + 1\}$,
and let $U_2$ be the set of vertices of degree $r-1$ in $H_2$.
Then $|U_2|= 2\rho - r + 2 = m(r - 1) + 1$.
In order to construct $G$,
we apply procedures (f) and (g) to $V(H_3) \setminus \{b_m\}$
while
\begin{itemize}
\item
joining $b_m$ to
$\lceil \frac{r-1}2 \rceil-1$ vertices of $U_1$ and
$\lfloor \frac{r-1}2 \rfloor + 1$ vertices of $U_2$.
\end{itemize}
We then obtain an $r$-regular graph $G'$
by joining each vertex of $U_1 \cup U_2$ to exactly one vertex of $H_3$.
Similarly to the argument in the case where both $m$ and $r$ are even,
it follows that $M$ does not extend to a perfect matching,
and the number of edge-disjoint odd ears of $V(M)$ is 
$2(\rho - \bup)+1 + m = (m+1)(r-1) - 2 \bup + m + 1 > mr-\bup$.
Since $\lambda_c(G') = |E(H_1, H_3)| = m(r - 1) - 1$,
the condition on $\lambda_c(G)$ in (i) of Theorem \ref{main} is sharp
in the case where either $m$ or $r$ is odd.

\subsection{Sharpness of the condition on $k$ in (ii)}
Let $m$ and $r$ be integers with $m \ge 2$ and $r \ge 3$,
and let $d$ be an integer with $d \ge (2m-1)(r-1)$.
Let $A$ be a graph isomorphic to $\Gamma(m, d; r)$, and
let $\{x_1, \ldots , x_m\}$ and $\{y_1, \ldots , y_m\}$ be the sets of vertices of degree $r-1$
in the two partite sets of $A$, respectively.
Let $G = A + \{x_1 x_2, y_1 y_2\} + \{x_i y_i \mid 3 \le i \le m\}$.
Since $\{x_1, \ldots , x_m\} \cup \{y_1, \ldots , y_m\}$ is
a good distance $d+1$ set of $A$ and $\lambda_c(A) \ge d$,
we can deduce that $\lambda_c(G) \ge d$.
Let $M= \{x_1 x_2, y_2 x_2^*\} \cup \{x_i y_i \mid 3 \le i \le m\}$,
where $x_2^*$ is a neighbor of $y_2$ in $G$ that is not $y_1$.
Since $G-V(M)$ is an unbalanced bipartite graph,
$M$ does not extend to a perfect matching in $G$.
\medskip

In what follows, we show that there exist
$(m-1)r$ edge-disjoint odd ears of $V(M)$ in $G$.
Let $M' = \{x_i y_i \mid 1 \le i \le m\}$ and $G' = A + M'$,
that is, $G' = G - \{x_1 x_2, y_1 y_2\} + \{x_1 y_1, x_2 y_2\}$.
Since $V(M')$ is a good distance $d+1$ set of $A$ and $\lambda_c(A) \ge d$,
we obtain $\lambda_c(G') \ge d \ge mr$.
Since $G'$ is an $r$-regular bipartite graph and
$M' \setminus \{x_1 y_1\}$ is a distance $d$ matching,
it follows from Lemma \ref{ears} that
there exists a set $\P_1$ of
$mr$ edge-disjoint odd ears of $V(M')$ in $G'$.
Since $G'$ is $r$-regular, every edge incident with $y_i$ is contained in
some path of $\P_1$ for each $i$.
Hence, there exists an ear $P \in \P_1$ containing the edge $y_2 x_2^*$.
By replacing $P$ with the path $y_2 x_2^*$ (of length one),
and removing some edges of the paths in $\P_1$ if necessary,
we obtain a set $\P_2$ of $mr$ edge-disjoint odd ears of
$V(M') \cup \{x_2^*\}$ in $G'$.
Then, by deleting from $\P_2$
the ear $x_2 y_2$ (of length one)
and all the ears containing the vertex $y_1$,
and adding the path $x_1 x_2$ (of length one),
we obtain $(m-1)r$ edge-disjoint odd ears of $V(M)$ in $G$.
Since $G$ is an $r$-regular graph with $\lambda_{oc}(G) \ge \lambda_c(G) \ge d$,
it follows that the condition $k=mr-r+1$ in (ii) of Theorem \ref{main} is sharp.
\section{Conclusion}
\label{conclusion}
In this paper, we have investigated the extendability of large matchings
in regular non-bipartite graphs,
which had previously been established only for bipartite graphs.
In particular, we have clarified that edge-disjoint odd ears
of the vertices of the matching to extend
play an important role in this setting,
and have obtained a sharp result 
with respect to the number of edge-disjoint odd ears (Theorem \ref{main} (i)).
Furthermore, by introducing a condition on odd-cyclic edge cuts,
we have obtained a result (Theorem \ref{main} (ii))
which extends the results of Plummer \cite{Pl91} and Aldred et al.~\cite{AJP}
for bipartite graphs to non-bipartite graphs, 
except for certain exceptional cases.
\medskip

Taking into account that no condition on $\lambda_c(G)$ is imposed in Theorem \ref{main} (ii), 
this result together with Lemma \ref{ears} provides a new perspective on
Theorems \ref{Pl91} and \ref{AJP}:
the condition on $\lambda_c(G)$ in Theorems \ref{Pl91} and \ref{AJP}
serves to ensure the existence of many edge-disjoint odd ears of $V(M)$,
while the assumption that $G$ is bipartite
is a condition to guarantee that $\lambda_{oc}(G)$ is sufficiently large.

\end{document}